# *On the Solutions of* $\left|\begin{array}{c}A\\X\end{array}\right| = \pm d$

*Massimo Salvi*

*Department of Mathematics and Computer Science*

*90123-Palermo, PA.*

*Italy*

e-mail: *massimo.salvi1@istruzione.it*

## *Abstract*

In this paper we deal with a non-linear Diophantine equation which arises from the determinant computation of an integer matrix. We show how to find a solution, when it exists. We define an equivalence relation and show how the set of all the solutions can be partitioned in a finite set of equivalence classes and find a set of solutions, one for each of these classes. We find a formula to express all the solutions and a formula to compute the cardinality of the set of fundamental solutions. An algorithm to compute the solutions is proposed and clarified with some examples.

## 1. Introduction

In this paper we deal with the diophantine equation (which has been already introduced in [1]):

$$\left|\begin{array}{c}A\\X\end{array}\right| = \pm d \tag{1}$$

where *A* and *X* are the matrices defined as follows:

*A* is a matrix $r \times n$ $\begin{bmatrix} a_{11} & a_{12} & .. & a_{1n} \\ .. & & & \\ a_{r1} & a_{r2} & .. & a_{rn} \end{bmatrix}$ with n columns and *r* rows, and $a_{ij} \in \mathbb{Z}$, similarly

*X* is a matrix $\begin{bmatrix} x_{r+1\,1} & x_{r+1\,2} & .. & x_{r+1\,n} \\ .. & & & \\ x_{n1} & x_{n2} & .. & x_{nn} \end{bmatrix}$ with n columns and *n-r* rows, and $x_{ij} \in \mathbb{Z}$.

then we denote by $\begin{bmatrix} A \\ X \end{bmatrix}$ the $n \times n$ matrix $\begin{bmatrix} a_{11} & a_{12} & .. & a_{1n} \\ .. & & & \\ a_{r1} & a_{r2} & .. & a_{rn} \\ x_{r+11} & x_{r+12} & .. & x_{r+1n} \\ .. & & & \\ x_{n1} & x_{n2} & .. & x_{nn} \end{bmatrix}$

Given any $n \times n$ matrix $M$, in this paper we denote by $|M|$ the determinant of $M$.

The equation (1) generalizes an equation studied by Smith [3] and, before him, by Hermite [8] in a more particular case. Smith analyzed the case in which $d=1$ and the greatest divisor of $A$ (see below) is 1, Hermite analyzed the case in which $d=1$ and $A$ is a vector of relatively prime integer numbers.

Recently, some algorithms have been developed in order to solve efficiently some diophantine problems, such as systems of linear diophantine equations arising from linear programming problems. Some of these algorithms use the Hermite Normal Form of integer matrices and the euclidean reduction algorithm applied to a vector of integers [2, 5, 6, 7, 13]. In this paper we use and develop some of this ideas to shed some light on the previous equation (1). In this paper we are not interested in the efficiency of the computation but only at some theoretical facts. First we clear up some notations that we use in the paper. We are interested in the entries $x_{ij}$ which satisfy $\left|\begin{matrix} A \\ X \end{matrix}\right| = d$ or $\left|\begin{matrix} A \\ X \end{matrix}\right| = -d$, for this reason we use the notation $\pm d$ in (1).

In this paper we give two main results:

1. After having defined an equivalence relation for the solutions of (1), we show that the set S of all the solutions can be decomposed into a finite number of equivalence sets, which form a partition of S. We show how to find a solution for each of these sets.
2. We give a formula which permits to write all the solutions of (1) in a convenient way.

In Section 2 we recall some concepts, definitions and theorems that are used in the rest of the paper. In Section 3 we introduce the definition of solution of (1) and the notion of L-A-equivalence. In Section 4 we recall a necessary and sufficient condition for the solvability of equation (1). In Section 5, we use the concept of L-A-equivalence to build a partition of the set of all solutions of (1), and we give a formula to express the solutions.

## 2. Some Preliminar Definitions and Results

In this section we recall some basic (and well known) results which are used in the following sections.

2.1 *Elementary Operations and Unimodular Matrix*

We recall some definitions and facts which are usefull to understand what follows in the paper.

Given an integer matrix *M*, an *elementary row operation* [12, p.12] is:

- Multiply a row by -1;
- Replace a row $R_i$ by $R_i + nR_j$, with $i \neq j$ and $n \in Z$     ;     (2)
- Swap two rows.

The *elementary column operations* are defined in the same way by replacing the rows with the columns. Note, anyway, that the swapping of two rows (or two columns ), say $R_i$ and $R_j$, can be obtained by using only the first two operations by the sequence of operations: $R_i \rightarrow R_i + R_j$ ; $R_j \rightarrow R_i - R_j$ ; $R_i \rightarrow R_i - R_j$.

An integer matrix *U*, such that $|U| = \pm 1$, is called *unimodular matrix* [2]. It's easy to prove that the inverse of an unimodular matrix is unimodular, the product of two unimodular matrices is uni modular, the identity matrix I is clearly unimodular. That is, the set of unimodular $n \times n$ matrices, for a given *n*, provided with the matrices multiplication operation is a group, the unimodular group, usually denoted by GL(n, Z). It is known that each of the above operations (2) is equivalent respectively to:

- pre-multiplying *M* by an unimodular matrix *U* if we operate on the rows
- post-multiplying *M* by an unimodular matrix *U* if we operate on the columns

A finite sequence of operations (2) on the rows or, respectively, on the colums of *M*, can be represented by an unimodular matrix, given by the product of all the matrices representing the elementary operations.

2.2 *Hermite Normal form of a Matrix*

**Definition 2.2.1** Given a full rank integer matrix *M*, we will say it is in *Hermite Normal Form* (HFN in what follows), if:

- The matrix *M* is of the form $[N \quad O]$ with *N* a nonsingular square matrix , and *O* is a matrix with all entries equal to 0;

- *N* is lower triangular ;

- All the entries of *N* are not negative and for every row of *N* the unique maximum is the diagonal entry which is greater than 0.

**Example 2.2.1**: $\begin{bmatrix} 1 & 0 & 0 & 0 & 0 \\ 0 & 3 & 0 & 0 & 0 \\ 1 & 1 & 2 & 0 & 0 \end{bmatrix}$

If *M* is a full rank matrix (let's consider the number of columns greater or equal to the number of rows), then we can define two different forms of HNF:

1. A first form defined as in Definition 2.2.1, which we refer to with notation **RHNF** (because, as claimed in the following Theorem 2.2.1, this form can be obtained by multiplying M on the right by an unimodular matrix ) .

   **Example 2.2.2**: $\begin{bmatrix} 1 & 0 & 0 \\ 1 & 2 & 0 \\ 1 & 2 & 3 \end{bmatrix}$

2. A second form defined as in Definition 2.2.1 but where for every column of N the unique maximum is the diagonal entry which is greater than 0. We refer to this form with notation **LHNF** (because, as we claim in the following Theorem 2.2.2, this form can be obtained by multiplying M on the left by an unimodular matrix).

   **Example 2.2.3**: $\begin{bmatrix} 1 & 0 & 0 \\ 0 & 3 & 0 \\ 0 & 2 & 2 \end{bmatrix}$

Some of the following theorems can be proved by applying the algorithms which we explain in 2.3. (see also [3, 8, 9]):

**Theorem 2.2.1** If $M$ is a full rank integer matrix with r rows and c columns and $r \leq c$, then there exists an unimodular matrix $U$ such that the matrix $MU$ is RHNF.

**Theorem 2.2.2** If $M$ is an integer square matrix, not singular, then there exists an unimodular matrix $U$ such that the matrix $UM$ is LHNF.

**Definition 2.2.2** If two integer square matrices $M$ and $M'$ are given and there exists an unimodular matrix $U$ such that $M=UM'$, then we say that $M$ and $M'$ are L-equivalent (left-equivalent).

**Definition 2.2.3** If two integer square matrices $M$ and $M'$ are given and there exists an unimodular matrix $U$ such that $M=UM'$, then we say that $M$ and $M'$ are R-equivalent (right-equivalent).

It's easy to prove that the previous definitions provide two equivalence relations.

The following theorem can be proved by multiplying and using the definitions [3]:

**Theorem 2.2.3** If $M$ and $M'$ are two integer square matrices, both not singular and reduced to RHNF, and R-equivalent, then $M=M'$.

In a similar way can be proved the following:

**Theorem 2.2.4** If *M* and *M'* are two integer square matrices, both not singular and reduced to LHNF, and L-equivalent, then *M=M'*.

2.3 *Greatest Divisor of a Matrix*

We can refer to Smith [3] to give the followings:

**Definition 2.3.1** If *M* is a matrix, the *determinants of M* are defined as the determinants of the greatest square matrices contained in *M*

**Definition 2.3.2** If *M* is an integer matrix, the *greatest divisor* of *M* is defined as the greatest common divisor of the determinants of *M*

Note that if *M* is a full rank matrix, then its greatest divisor is well defined and different from 0. If *M* is not full rank, then all the determinants are 0, so the greatest divisor is not defined.

**Definition 2.3.3** If *M* is an integer matrix, then *M* is said *prime matrix* if its greatest divisor is 1.

**Theorem 2.3.1** If *M* is a full rank integer matrix, then, by pre-or post-multiplying *M* by an unimodular matrix, the greatest divisor of the resulting matrix is the same as the one of *M*.

See [3] for a proof.

If *M* is in HNF (RHNF or LHNF), that is $M=[N \quad O]$, it's clear that the greatest divisor of *M* can be obtained by the product of all the diagonal entries of the square matrix *N*. In fact this is the only square matrix contained in *M* with determinant not equal to 0. Hence, by using the previous theorems 2.4.1 and 2.2.1 we have the following:

**Theorem 2.3.2** If *M* is a full rank integer matrix with r rows and c columns and $r \leq c$, and *M'*= $[N \quad O]$ is the HNF of *M*, then the greatest divisor of *M* is equal to the greatest divisor of *M'*, which is equal to the product of all the diagonal entries of *N*.

2.4 *Algorithm to reduce a Matrix in HNF*

We show an algorithm which permits to reduce a matrix in RHNF by using the elementary operations (2) on the columns. Similar algorithms have been presented by various authors, [2, 4, 5, 6] differing from each other by the different complexity of computation. In this paper we are not interested in the efficiency of the computation but we only want to show a way to do it. The algorithm is based on the Extended Euclidean Algorithm (see [7, 13]) which we can use to obtain an integer column (or row) reduction in a integer matrix *M*. We explain the algorithm and use it to reduce to Right Hermite Normal Form a full rank integer matrix *M* with *r* rows and *c* columns, with columns and $r \leq c$.

Be $M = \begin{bmatrix} m_{11} & m_{12} & .. & .. & .. & m_{1c} \\ m_{21} & .. & .. & .. & .. & .. \\ .. & .. & .. & .. & .. & .. \\ m_{r1} & .. & .. & .. & .. & m_{rc} \end{bmatrix}$

If we consider the first row, all its entries cannot be all equal to 0 since the *M* is full rank, so the greatest common divisor $d_1$ of its *c* elements $m_{1i}$ exists. We can reduce *M* into the following form *M'* by unimodular operations on the columns.

$$M' = \begin{bmatrix} d_1 & 0 & .. & .. & .. & 0 \\ m'_{21} & .. & .. & .. & .. & .. \\ .. & .. & .. & .. & .. & .. \\ m'_{r1} & .. & .. & .. & .. & m'_{rc} \end{bmatrix}$$

Algorithm 1

To obtain this form we can repeat the following algorithm until in the first row there is only one non-zero entry, $d_1$:

1. Find the value of the row that has the smaller absolute value different from zero, and we denote it by $m_{1j}$

2. For every $i \neq j$ using the division algorithm we can write:

$m_{1i} = q_i m_{1j} + r_i$ with $0 \leq r_i < |m_{1j}|$

3. In each column $C_i$, $i \neq j$, we operate the substitution: $C_i \to C_i - q_i C_j$

it can be proved that $d_1$ is the greatest common divisor of the entries of the first row (see [7]). It can be seen as a particular case of theorem 2.3.1 where the matrix has only one row. We can move $d_1$ in the first position of the row, eventually by swapping two columns or changing the signs of one, to obtain a matrix in the form *M'*. Then we can continue this process starting from the matrix *M'*, without considering the first column and considering the entries of the second row. By this way we will be able to transform the matrix *M'* in a matrix *M''* of the form:

$$M'' = \begin{bmatrix} d_1 & 0 & .. & .. & .. & 0 \\ m'_{21} & d_2 & 0 & .. & .. & 0 \\ .. & .. & .. & .. & .. & .. \\ m'_{r1} & m''_{r2} & .. & .. & .. & m''_{rc} \end{bmatrix}$$

Going on, in this way we can transform the matrix into a lower triangular matrix. All the operations are equivalent to post-multiply *M* by an unimodular matrix *U*. Note that for all the diagonal entries we have $d_i > 0$, because if $d_i < 0$ we can multiply the column by -1, and for every *i*, $d_i \neq 0$ otherwise, for theorem 2.3.1 and theorem 2.3.2, the greatest divisor of *M* cannot be defined and positive as required by the hypotesis that *M* is full rank.

Algorithm 2

In order to transform the matrix in the form *M''* in Hermite Normal Form we can apply the following algorithm :

1. We start from *i*=2 ;

2. For every $j < i$, by using the division algorithm we can write:

$m_{ij} = q_j d_i + R_{ij}$ with $0 \leq R_{ij} < d_i$;

3. In each column $C_j$, $j < i$, we operate the substitution: $C_j \to C_j - q_j C_i$;

4. Increase $i$ and repeat form point 1 untill $i = r$.

At the end of algorithm 2 the matrix will be in the Right Hermite Normal Form, that is:

$$M''' = \begin{bmatrix} d_1 & 0 & .. & .. & .. & 0 \\ R_{21} & d_2 & 0 & .. & .. & 0 \\ .. & .. & .. & .. & .. & .. \\ R_{r1} & R_{r2} & .. & d_r & .. & 0 \end{bmatrix}$$

We can obviously modify these algorithms to operate on rows instead that on the columns, in order to place the full rank matrix $M$ in LHNF. In this way we get a proof of theorems 2.2.1 and 2.2.2.

Example 2.4.1

Let's consider the full rank matrix $\begin{bmatrix} 2 & 2 & -3 & 4 \\ 2 & 2 & 1 & 2 \end{bmatrix}$ and apply the previous algorithms.

| | Column transformation | Equivalent matrix post-multiplication |
|---|---|---|
| | Algorithm 1 | |
| T1 | $\begin{cases} C2 \to C2 - C1 \\ C3 \to C3 + C1 \\ C4 \to C4 - 2C1 \end{cases}$ | $\begin{bmatrix} 2 & 2 & -3 & 4 \\ 2 & 2 & 1 & 2 \end{bmatrix} \cdot \begin{bmatrix} 1 & -1 & 1 & 2 \\ 0 & 1 & 0 & 0 \\ 0 & 0 & 1 & 0 \\ 0 & 0 & 0 & 1 \end{bmatrix} = \begin{bmatrix} 2 & 0 & -1 & 0 \\ 2 & 0 & 3 & -2 \end{bmatrix}$ |
| T2 | $\{C3 \to -C3$ | $\begin{bmatrix} 2 & 0 & -1 & 0 \\ 2 & 0 & 3 & -2 \end{bmatrix} \cdot \begin{bmatrix} 1 & 0 & 0 & 0 \\ 0 & 1 & 0 & 0 \\ 0 & 0 & -1 & 0 \\ 0 & 0 & 0 & 1 \end{bmatrix} = \begin{bmatrix} 2 & 0 & 1 & 0 \\ 2 & 0 & -3 & -2 \end{bmatrix}$ |
| T3 | $\{C1 \to C1 - 2C3$ | $\begin{bmatrix} 2 & 0 & 1 & 0 \\ 2 & 0 & -3 & -2 \end{bmatrix} \cdot \begin{bmatrix} 1 & 0 & 0 & 0 \\ 0 & 1 & 0 & 0 \\ -2 & 0 & -1 & 0 \\ 0 & 0 & 0 & 1 \end{bmatrix} = \begin{bmatrix} 0 & 0 & 1 & 0 \\ 8 & 0 & -3 & -2 \end{bmatrix}$ |
| T4 | $\{C1 \leftrightarrow C3$ | $\begin{bmatrix} 0 & 0 & 1 & 0 \\ 8 & 0 & -3 & -2 \end{bmatrix} \cdot \begin{bmatrix} 0 & 0 & 1 & 0 \\ 0 & 1 & 0 & 0 \\ 1 & 0 & 0 & 0 \\ 0 & 0 & 0 & 1 \end{bmatrix} = \begin{bmatrix} 1 & 0 & 0 & 0 \\ -3 & 0 & 8 & -2 \end{bmatrix}$ |

| T5 | $\{C3 \rightarrow C3+4C4$ | $\begin{bmatrix} 1 & 0 & 0 & 0 \\ -3 & 0 & 8 & -2 \end{bmatrix} \cdot \begin{bmatrix} 1 & 0 & 0 & 0 \\ 0 & 1 & 0 & 0 \\ 0 & 0 & 1 & 0 \\ 0 & 0 & 4 & 1 \end{bmatrix} = \begin{bmatrix} 1 & 0 & 0 & 0 \\ -3 & 0 & 0 & -2 \end{bmatrix}$ |
| --- | --- | --- |
| T6 | $\begin{cases} C2 \rightarrow -C4 \\ C4 \rightarrow C2 \end{cases}$ | $\begin{bmatrix} 1 & 0 & 0 & 0 \\ -3 & 0 & 0 & -2 \end{bmatrix} \cdot \begin{bmatrix} 1 & 0 & 0 & 0 \\ 0 & 0 & 0 & 1 \\ 0 & 0 & 1 & 0 \\ 0 & -1 & 0 & 0 \end{bmatrix} = \begin{bmatrix} 1 & 0 & 0 & 0 \\ -3 & 2 & 0 & 0 \end{bmatrix}$ |
| | Algorithm 2 | |
| T7 | $\{C1 \rightarrow C1+2C2$ | $\begin{bmatrix} 1 & 0 & 0 & 0 \\ -3 & 2 & 0 & 0 \end{bmatrix} \cdot \begin{bmatrix} 1 & 0 & 0 & 0 \\ 2 & 1 & 0 & 0 \\ 0 & 0 & 1 & 0 \\ 0 & 0 & 0 & 1 \end{bmatrix} = \begin{bmatrix} 1 & 0 & 0 & 0 \\ 1 & 2 & 0 & 0 \end{bmatrix}$ |

Now, by multiplying all the unimodular matrices we obtain :

$$\begin{bmatrix} 2 & 2 & -3 & 4 \\ 2 & 2 & 1 & 2 \end{bmatrix} \cdot \begin{bmatrix} 3 & 2 & -5 & -1 \\ 0 & 0 & 0 & 1 \\ -1 & 0 & 2 & 0 \\ -2 & -1 & 4 & 0 \end{bmatrix} = \begin{bmatrix} 1 & 0 & 0 & 0 \\ 1 & 2 & 0 & 0 \end{bmatrix} \qquad (3)$$

### 3. Equivalent Solutions of (1)

In this section we introduce an equivalence relation between two solutions of (1), which we will use in the following part of the paper. First we define a solution of (1):

**Definition 3.1** Let's suppose that a set of integer entries $x_{ij}$ satisfy the equation $\begin{vmatrix} A \\ X \end{vmatrix} = d$ or $\begin{vmatrix} A \\ X \end{vmatrix} = -d$, then, by denoting with $X$ the matrix of the entries $x_{ij}$, that is $X = \begin{bmatrix} x_{k+1\,1} & x_{k+1\,2} & .. & x_{k+1\,n} \\ .. & & & \\ x_{n1} & x_{n2} & .. & x_{nn} \end{bmatrix}$, we will define $X$ *a solution* of (1).

**Definition 3.2** If $S_i$ and $S_j$ are two solutions (we suppose that in the equation (1) the matrix $A$ is given), we say that $S_i$ is *L-A- equivalent to* $S_j$, if there exists an unimodular matrix $U$ such that:

$$\begin{bmatrix} A \\ S_i \end{bmatrix} = U \begin{bmatrix} A \\ S_j \end{bmatrix}$$

It can be easily proved that definition 3.2 provides an *equivalence relation*:

- Reflexivity: $\begin{bmatrix} A \\ S_i \end{bmatrix} = I \begin{bmatrix} A \\ S_i \end{bmatrix}$;

- Simmetry: if $U$ exists satisfying $\begin{bmatrix} A \\ S_i \end{bmatrix} = U \begin{bmatrix} A \\ S_j \end{bmatrix}$, then by multiplying by $U^{-1}$, which is unimodular, we obtain $U^{-1} \begin{bmatrix} A \\ S_i \end{bmatrix} = \begin{bmatrix} A \\ S_j \end{bmatrix}$;

- Transitivity: if $U$ exists satisfying $\begin{bmatrix} A \\ S_i \end{bmatrix} = U \begin{bmatrix} A \\ S_j \end{bmatrix}$, and $U'$ exists satisfying $\begin{bmatrix} A \\ S_j \end{bmatrix} = U' \begin{bmatrix} A \\ S_k \end{bmatrix}$, then $\begin{bmatrix} A \\ S_i \end{bmatrix} = UU' \begin{bmatrix} A \\ S_k \end{bmatrix}$, where $UU'$ is unimodular.

### 4. A necessary and sufficient Condition for Solvability of (1)

In [1] it has been proved that, by denoting with $|A|$ the greatest divisor of $A$, holds the following:

**Theorem 4.1** If $A$ is given, the equation 1.1 ($d > 0$) has integer solutions if and only if $A$ is full rank and $|A|$ divides $d$. If $A$ is not given the equation (which becames $|X| = \pm d$) is always solvable, even if $d=0$.

The proof of this theorem is implicitly embedded in the following proof of Theorem 5.2.2.

**Corollary 4.1** If tha matrix $A$ is prime, then the equation (1) is always solvable.

**Example 4.1**

Be the equation to solve: $\begin{vmatrix} 1 & 2 & -3 & 4 \\ 0 & 1 & 1 & 2 \\ x_{31} & x_{32} & x_{33} & x_{34} \\ x_{41} & x_{42} & x_{43} & x_{44} \end{vmatrix} = \pm 2$

The greatest divisor of the matrix $A$ is g.c.d $\{1,1,2,5,0,-10\} = 1$, namely $A$ is a prime matrix. This means, by theorem 4.1, that the equation is solvable.

**Example 4.2**

Be the equation to solve: $\begin{vmatrix} 2 & 2 & -3 & 4 \\ 2 & 2 & 1 & 2 \\ x_{31} & x_{32} & x_{33} & x_{34} \\ x_{41} & x_{42} & x_{43} & x_{44} \end{vmatrix} = \pm 3$

The greatest divisor of the matrix $A$ is g.c.d $\{0,8,-4,8,-4,-10\} = 2$. Clearly 2 does not divides 3, so, by theorem 4.1, the equation is not solvable.

# 5. The Solutions of $\begin{vmatrix} A \\ X \end{vmatrix} = \pm d$

We start by considering a simpler case which permits us to introduce the way we will operate to solve the general case.

## 5.1 The case $|X| = \pm d$

First let's consider the simpler case in which there are no constants, namely the matrix $A$ is not given.

**Theorem 5.1.1** Given the integers $d > 0$ and $n > 0$, there is only a finite number of $n \times n$ integer matrix in LHNF (or RHNF) with determinant equal to $d$.

*Proof.* Be $S_i$ a matrix in LHNF such that $|S_i| = d$. The matrix $S_i$ is lower triangular and in the form: $\begin{bmatrix} a_{11} & 0 & 0 & .. & 0 \\ .. & .. & 0 & .. & 0 \\ a_{ij} & & a_{ii} & .. & 0 \end{bmatrix}$ where $a_{ii} > 0 \; \forall i$ and $0 \le a_{ij} < a_{jj} \; \forall i, \forall j$ (4)

from $|S_i| = d$ we have that $\prod_{j=1}^{n} a_{jj} = d$ (5)

but since $a_{jj}$ and $n$ are integers, the previous equation (5) can be satisfyed only by a finite number of values for $a_{jj}$. From condition we have $a_{ii} > 0 \; \forall i$ and $0 \le a_{ij} < a_{jj} \; \forall i, \forall j$, therefore for every value of $a_{jj}$, only a finite number of integer values is allowed for each $a_{ij}$. We can conclude that only a finite number of matrix of the form (4) exists with determinant equal to $d$. Similarly we can prove the same results also for matrices in RHNF. □

**Definition 5.1.1** If $d > 0$ and $n > 0$ are integers, and $S_i$ is a $n \times n$ matrix in LHNF such that $|S_i| = d$ we define the *class of L-equivalence* of $S_i$ as the set $C_i$ of all the matrices which are L-equivalent to $S_i$.

From theorem 5.1.1 it follows that, given the integers $d > 0$ and $n > 0$, there is a finite number of classes $C_i$ of L-equivalence.

**Theorem 5.1.2** If $d > 0$ and $n > 0$ are integers, and we denote with $S$ the set of the all the solutions of the equation $|X| = \pm d$, where $X$ is a $n \times n$ integer matrix, then the set of classes of L-equivalence defined as in definition 5.1.1 is a partition of $S$.

*Proof.* Given a solution $X$, for theorem 2.2.2 there exists an unimodular matrix $U$ such that $UX = S_i$, where $S_i$ is LHNF and $|S_i| = d$, namely $X$ is L-equivalent to $S_i$. This means that, by definition 5.1.1, $X$ belongs to $C_i$.

If we suppose that $X$ belongs to two differents classes, say $C_i$ and $C_j$, then there exist two unimodular matrices $U_i$ and $U_j$ such that $U_i X = S_i$ and $U_j X = S_j$. Hence we can write $X = U_i^{-1} S_i$ and $X = U_j^{-1} S_j$, and therefore

$$U_i^{-1} S_i = U_j^{-1} S_j$$

$$S_i = U_i U_j^{-1} S_j$$

The matrix $U_i U_j^{-1}$ is unimodular, hence $S_i$ is L-equivalent to $S_j$, but for theorem 2.2.4 we have $S_i = S_j$ and consequently $C_i = C_j$. □

5.2 *The general case* $\left| \begin{array}{c} A \\ X \end{array} \right| = \pm d$

We generalize the previous theorem 5.1.2 by showing how to build a finite number of solutions $S_i$ of equation (1), which play the same role as in the simpler case $|X| = \pm d$ analyzed in the previous section.

5.2.1 *Construction of a Set of Solutions*

Let's consider the equation $\left| \begin{array}{c} A \\ X \end{array} \right| = \pm d$, and suppose it's solvable. For theorem 4.1 we have that $|A|$ is defined, not equal to 0, it divides $d$, and $A$ is full rank. For theorem 2.2.1 there exists an unimodular matrix $U$ such that $AU$ is RHNF. We can write more explicitly such form as follows:

$$AU = \begin{bmatrix} a_{11} & 0 & 0 & .. & 0 \\ .. & .. & 0 & .. & 0 \\ a_{ij} & & a_{ii} & .. & 0 \end{bmatrix} \text{ where } a_{ii} > 0 \; \forall i \text{ and } 0 \leq a_{ij} < a_{ii} \; \forall i, \forall j \; , j < i$$

(7)

In order to proceed in the construction we need the following:

**Theorem 5.2.1** If a full rank matrix $A$ is in the form $[N \quad O]$, where $N$ is a square lower triangular matrix and $O$ is a matrix with all entries equal to 0, then it is possible to find an unimodular matrix $U'$ such that $U'A$ is in the form :

$$[N' \ O] = \begin{bmatrix} a'_{11} & 0 & 0 & .. & 0 \\ .. & .. & 0 & .. & 0 \\ a'_{ij} & & a'_{ii} & .. & 0 \end{bmatrix} \text{ where } a'_{ii} > 0 \ \forall i \text{ and } 0 \leq a'_{ji} < a'_{ii} \ \forall i, \forall j, j < i \qquad (8)$$

(That is, in (8) for each column the maximum is the diagonal entry, whereas in (7) it is true for each row).

*Proof.* The square matrix *N'* in (8) is not singular, hence we can apply theorem 2.2.2 and find an unimodular matrix *U'* such that *U' N'* is in LHNF. Then, if we consider the product, we obtain:

$$U'A = [U'N' \ U'O] = [U'N' \ O] = \begin{bmatrix} a'_{11} & 0 & 0 & .. & 0 \\ .. & .. & 0 & .. & 0 \\ a'_{ij} & & a'_{ii} & .. & 0 \end{bmatrix}$$

and the conditions (8) are satisfyed since the matrix *U' N'* is in LHNF. □

Now we take into consideration the matrix *U'AU*, where *U* and *U'* are defined as in (7) and (8):

$$U'AU = \begin{bmatrix} a'_{11} & 0 & 0 & .. & 0 \\ .. & .. & 0 & .. & 0 \\ a'_{i1} & & a'_{ii} & .. & 0 \end{bmatrix} \qquad (9)$$

From theorems 2.3.1 and 2.3.2 we have that the greatest divisor of the matrix in (9),

is equal to the greatest divisor of *A*:

$$\prod_i a'_{ii} = |U'AU| = |A| \qquad (10)$$

from (10) it follows that $\prod_i a'_{ii}$ divides *d*. Therefore we can consider the integer *k* :

$$k = \frac{d}{|A|} \qquad (11)$$

Now let's consider the set of $n \times n$ matrices in LHNF in which $\prod_j b_{jj} = k$ : \qquad (12)

$$\begin{bmatrix} U'AU \\ B_i \end{bmatrix} = \begin{bmatrix} a'_{11} & 0 & .. & .. & .. & 0 \\ .. & .. & 0 & .. & .. & 0 \\ a'_{r1} & & a'_{rr} & 0 & .. & 0 \\ b_{l1} & b_{l2} & .. & b_{ll} & .. & 0 \\ .. & .. & & & .. & 0 \\ b_{n1} & & & & & b_{nn} \end{bmatrix} \qquad (13)$$

Note that, by construction, the determinant of such matrices is $\prod_{i=1}^{r} a_{ii} \cdot \prod_{j=r+1}^{n} b_{jj} = d$, where $r$ and $n$ are, respectively, the number of row and of column of $U'AU$ (the same as $A$).

Note also that, from theorem 5.2.1, the matrix $U'AU$ is already in the correct form to give rise to a LHNF.

Since we suppose the matrices in (13) be in LHNF, the entries of $B_i$ have to satisfy :

$$0 \leq b_{ij} < a'_{jj} \quad \forall i \; \forall j, \; j < i \quad \text{if} \quad j \leq r$$

$$0 \leq b_{ij} < b_{jj} \quad \forall i \; \forall j, \; j < i \quad \text{if} \quad j > r \tag{14}$$

From the previous conditions (12) and (14) it follows that there exists only a finite number of matrices $B_i$ (see theorem 5.1.1).

The next step is to build, starting from the matrices (13), a set of solutions of (1).

**Theorem 5.2.2** For each matrix $B_i$, defined as in (12), (13) and (14), and given the unimodular matrix $U^{-1}$, where $U$ is defined as in (7), the matrix $B_i U^{-1}$ is a solution of the equation $\begin{vmatrix} A \\ X \end{vmatrix} = \pm d$.

*Proof.* Starting from the unimodular $r \times r$ matrix $U'$, as in (9), let's define an unimodular $n \times n$ matrix $U''$ as follows:

$$U'' = \begin{bmatrix} U'^{-1} & O \\ O & I \end{bmatrix} \tag{15}$$

where $I$ is an unit $(n-r) \times (n-r)$ matrix, and $O$ are matrices with all entries equal to 0.

Now let's consider the product:

$$U'' \cdot \begin{bmatrix} U'AU \\ B_i \end{bmatrix} \cdot U^{-1} = \tag{16}$$

$$\begin{bmatrix} U'^{-1} & O \\ O & I \end{bmatrix} \begin{bmatrix} U'AU \\ B_i \end{bmatrix} \cdot U^{-1} = \begin{bmatrix} U'^{-1} & O \\ O & I \end{bmatrix} \begin{bmatrix} U'AUU^{-1} \\ B_i U^{-1} \end{bmatrix} = \begin{bmatrix} U'^{-1} & O \\ O & I \end{bmatrix} \begin{bmatrix} U'A \\ B_i U^{-1} \end{bmatrix} \tag{17}$$

note that the dimensions of the inner blocks of the matrices allow a block multiplication, hence the expression (17) is equal to:

$$= \begin{bmatrix} U'^{-1}U'A + OB_iU^{-1} \\ OU'A + IB_iU^{-1} \end{bmatrix} = \begin{bmatrix} A \\ B_iU^{-1} \end{bmatrix} \quad (18)$$

As we have shown in (13), the determinant of the matrix $\begin{bmatrix} U'AU \\ B_i \end{bmatrix}$ is $d$, the matrix defined in 5.2.9 is clearly unimodular and its determinant is 1 or -1, therefore the determinant of matrix (18) is $d$ or $-d$. We can conclude that the matrix $B_iU^{-1}$ is a solution of equation (1). □

Note that, from (16), the sign of the determinant of the matrix $\begin{bmatrix} U'AU \\ B_i \end{bmatrix}$ is the same for every possible $B_i$.

In what follows we define $S_i = B_iU^{-1}$ \quad (19)

We will refer to this set of solutions as *fundamental solutions* of the equation $\begin{vmatrix} A \\ X \end{vmatrix} = \pm d$

Now, similarly to what we have done in the simpler case previously analyzed (definition 5.1.1), we can give the following:

**Definition 5.2.3** If $S_i$ is a solution of $\begin{vmatrix} A \\ X \end{vmatrix} = \pm d$, $d > 0$, defined as in theorem 5.2.2, we define as *class of L-A-equivalence* of $S_i$ the set $C_i$ of all the matrices which are L-A-equivalent to $S_i$.

We have seen that there exists only a finite number of $B_i$, therefore there exists only a finite number of $S_i$ and of classes $C_i$.

5.2.2 *A Generalization of Theorem 2.2.4*

To proceed we need a generalization of theorem 2.2.4, which is given by the following:

**Theorem 5.2.3** If the square integer matrix $\begin{bmatrix} A \\ B \end{bmatrix}$ is L-equivalent to a square integer and not singular matrix $R$, which is in LHNF, and $A$ is a full rank matrix with $r$ rows in the form:

$$[N \quad O] = \begin{bmatrix} a_{11} & 0 & 0 & .. & 0 \\ .. & .. & 0 & .. & 0 \\ a_{r1} & & a_{rr} & .. & 0 \end{bmatrix}$$

in which $O$ is a matrix with all entries equal to 0 and $N$ is a square matrix in LHNF

then the matrix $A$ is equal to the matrix obtained by the first $r$ rows of the matrix $R$.

*Proof.* Let's consider the matrix $B$: it is full rank, because the matrix $\begin{bmatrix} A \\ B \end{bmatrix}$ is L-equivalent to $R$, hence its determinant is not 0. We first show, by applying on the rows of $B$ the algorithm yet explained in 2.4, that an unimodular matrix $H$ exists, such that $H$ is in the form:

$$H = \begin{bmatrix} I & O \\ L & H' \end{bmatrix} \quad (20)$$

and the following product is in LHNF

$$H \begin{bmatrix} A \\ B \end{bmatrix} = \begin{bmatrix} I & O \\ N & H' \end{bmatrix} \begin{bmatrix} A \\ B \end{bmatrix} = \begin{bmatrix} A \\ NA + H'B \end{bmatrix} \quad (21)$$

Given the matrix $B = \begin{bmatrix} b_{11} & .. & .. & .. & b_{1n} \\ .. & .. & .. & .. & .. \\ b_{k1} & .. & .. & .. & b_{kn} \end{bmatrix}$

the first step is to apply the algorithm 1 (section 2.4) to the first $k$ rows of $B$. By representing these operations by the unimodular matrix $U_1$, we can write:

$$\begin{bmatrix} b'_{11} & .. & .. & .. & 0 \\ .. & .. & .. & .. & 0 \\ b'_{k1} & .. & .. & .. & b'_{kn} \end{bmatrix} = U_1 B \quad (22)$$

then we can apply again the algorithm to the first $k$-1 rows of the matrix $U_1 B$. We can represent these operations by another unimodular matrix $U_2$, to obtain a new matrix in the form:

$$\begin{bmatrix} b''_{11} & .. & .. & 0 & 0 \\ .. & .. & .. & b_{k-1n-1} & 0 \\ b'_{k1} & .. & .. & .. & b'_{kn} \end{bmatrix} = U_2 U_1 B \quad (23)$$

so, going on in this way we obtain a matrix in the form

$$B' = [B_1 \quad B_2] = U_r .. U_2 U_1 B \quad (24)$$

in which $B_2$ is square lower triangular.

The second step is to apply the algorithm 2 on the rows of $B'$; we can represent these operations by the unimodular matrix $U'$, to obtain a matrix in the form:

$$B'' = U'B' = [B_1' \quad B_2'] = U' U_r .. U_2 U_1 B \quad (25)$$

in which the matrix $B_2'$ is in LHNF.

The third step is to consider the matrix $\begin{bmatrix} A \\ B'' \end{bmatrix}$ and apply the algorithm 2 on the rows of $B''$ and using the rows of $A$ to obtain, finally, a matrix in LHNF. Also in this case we can represent the algorithm by pre-multiplying by an unimodular matrix $U''$:

$$\begin{bmatrix} I & O \\ U''_1 & U''_2 \end{bmatrix} \begin{bmatrix} A \\ B'' \end{bmatrix} = \begin{bmatrix} A \\ U''_1 A + U''_2 B'' \end{bmatrix} \qquad (26)$$

If we denote the product in (25) by $U$, that is:

$$U = U' U_r .. U_2 U_1 \qquad (27)$$

we can synthetize all the operations with a unique uni modular matrix given by the product:

$$\begin{bmatrix} I & O \\ U''_1 & U''_2 \end{bmatrix} \begin{bmatrix} I & O \\ O & U \end{bmatrix} \begin{bmatrix} A \\ B \end{bmatrix} = \begin{bmatrix} I & O \\ U''_1 & U''_2 U \end{bmatrix} \begin{bmatrix} A \\ B \end{bmatrix} \qquad (28)$$

The matrix $H$ in (20) coincides with the one of the previous (28) after having set

$L = U''_1$ and $H' = U''_2 U$

Now, going back to the hypotesis of the theorem, an unimodular matrix $W$ exists such that:

$$\begin{bmatrix} A \\ B \end{bmatrix} = WR \qquad (29)$$

and $R$ is in LHNF.

From (21) we obtain:

$$\begin{bmatrix} A \\ B \end{bmatrix} = H^{-1} \begin{bmatrix} A \\ NA + H'B \end{bmatrix} \qquad (30)$$

by comparing (30) with (29) we can write:

$$WR = H^{-1} \begin{bmatrix} A \\ NA + H'B \end{bmatrix} \qquad (31)$$

and then

$$R = W^{-1} H^{-1} \begin{bmatrix} A \\ NA + H'B \end{bmatrix} \qquad (32)$$

In the previous (32) the matrices $R$ and $\begin{bmatrix} A \\ NA+H'B \end{bmatrix}$ are both in LHNF and the matrix $W^{-1}H^{-1}$ is unimodular, therefore, due to theorem 2.2.4, the two matrices are equal, in particular $A$ is equal to the first $r$ rows of the matrix $R$. □

### 5.2.3 A Partition of the Set of Solutions

**Theorem 5.2.4** Given the integers $d > 0$ and $n > 0$, and by denoting with $S$ the set of all the solutions of the equation (1), $\begin{vmatrix} A \\ X \end{vmatrix} = \pm d$, where $A$ is a full rank integer matrix, we have that the set of classes $C_i$ of L-A-equivalence defined as in definition 5.2.3 is a partition of $S$.

*Proof.* The first step is to prove that any solution of equation (1) is L-A-equivalent to one of the $S_i$.

If $X$ is a solution, that is :

$$\begin{vmatrix} A \\ X \end{vmatrix} = \pm d \tag{33}$$

by using the unimodular matrices $U$ and $U'$, defined as in (7) and theorem 5.2.1, we can consider the matrix:

$$\begin{bmatrix} U' & O \\ O & I \end{bmatrix} \cdot \begin{bmatrix} A \\ X \end{bmatrix} U = \begin{bmatrix} U'AU \\ XU \end{bmatrix} \tag{34}$$

in which the the determinant is $+d$ or $-d$ and the matrix $U'AU$ satisfies the conditions (8).

The matrix (34) is a solution of the equation , $|X'| = \pm d$ hence, by theorem 5.1.2, there exists an unimodular matrix $W$, such that:

$$\begin{bmatrix} U'AU \\ XU \end{bmatrix} = W R_i \tag{35}$$

where $R_i$ is LHNF and $|R_i| = d$

Since the matrix (34) satisfies the conditions (8), for theorem 5.2.3, we have that the matrix $R_i$, can be written as:

$$\begin{bmatrix} U'AU \\ R'_i \end{bmatrix} \tag{36}$$

in which the entries $r'_{ij}$ of the matrix $R'_i$ satisfy the conditions (similarly to (14)):

$0 \leq r'_{ij} < a'_{jj} \quad \forall i \text{ if } j \leq row$

$$0 \leq r'_{ij} < r'_{jj} \quad \forall i \quad \text{if} \quad j > row \tag{37}$$

in which we denoted by "*row*" the number of rows, and $a'_{jj}$ are the diagonal entries of the matrix $U'AU$.

Now, for what we have proved in theorem 5.2.2, in particular from (16) to (18), the matrix $\begin{bmatrix} U'AU \\ B \end{bmatrix}$, introduced in (13) is exactly one of the $R_i$ (is LHNF and its determinant is *d*); it means

that we can find a $B_i$ such that $R_i = \begin{bmatrix} U'AU \\ B_i \end{bmatrix}$ \hfill (38)

therefore, by using theorem 5.2.2 and by defining $U''$ as in (15), there exists a solution $S_i$ of $\begin{vmatrix} A \\ X \end{vmatrix} = \pm d$ such that:

$$U'' \cdot \begin{bmatrix} U'AU \\ B_i \end{bmatrix} \cdot U^{-1} = U'' R_i U^{-1} = \begin{bmatrix} A \\ S_i \end{bmatrix} \quad \text{(for (38))} \tag{39}$$

hence from (39) we can write:

$$R_i = U''^{-1} \begin{bmatrix} A \\ S_i \end{bmatrix} U \tag{40}$$

by substituting in (35) we obtain:

$$\begin{bmatrix} U'AU \\ XU \end{bmatrix} = W R_i = W U''^{-1} \begin{bmatrix} A \\ S_i \end{bmatrix} U \tag{41}$$

and post-multiplying by $U^{-1}$

$$\begin{bmatrix} U'A \\ X \end{bmatrix} = W R_i = W U''^{-1} \begin{bmatrix} A \\ S_i \end{bmatrix} \tag{42}$$

the first member of (42) can be written, by using the definition 5.2.9:

$$\begin{bmatrix} U'A \\ X \end{bmatrix} = \begin{bmatrix} U' & O \\ O & I \end{bmatrix} \cdot \begin{bmatrix} A \\ X \end{bmatrix} = U''^{-1} \cdot \begin{bmatrix} A \\ X \end{bmatrix} \tag{43}$$

so we obtain:

$$U''^{-1} \cdot \begin{bmatrix} A \\ X \end{bmatrix} = W U''^{-1} \begin{bmatrix} A \\ S_i \end{bmatrix} \tag{44}$$

and then:

$$\begin{bmatrix} A \\ X \end{bmatrix} = U''WU''^{-1} \begin{bmatrix} A \\ S_i \end{bmatrix} \tag{45}$$

in which the matrix $U''WU''^{-1}$ is unimodular.

The previous (45) means exactly that the solution $X$ is L-$A$-equivalent to one of the fundamental solutions $S_i$ defined in theorem 5.2.2.

The second step is to prove that a solution $X$ cannot be L-$A$-equivalent to two differents fundamental solutions, say $S_i$ and $S_j$.

If a solution $X$ is L-$A$-equivalent to $S_i$ and $S_j$, we have that:

$$\begin{bmatrix} A \\ X \end{bmatrix} = W \begin{bmatrix} A \\ S_i \end{bmatrix} \tag{46}$$

$$\begin{bmatrix} A \\ X \end{bmatrix} = W' \begin{bmatrix} A \\ S_j \end{bmatrix} \tag{47}$$

and from the previous (46) and (47) it follows:

$$W \begin{bmatrix} A \\ S_i \end{bmatrix} = W' \begin{bmatrix} A \\ S_j \end{bmatrix} \tag{48}$$

$$\begin{bmatrix} A \\ S_i \end{bmatrix} = W^{-1}W' \begin{bmatrix} A \\ S_j \end{bmatrix} \tag{49}$$

namely, $S_i$ is L-$A$-equivalent to $S_j$. From (19), it means that $B_i$ and $B_j$ exist, such that:

$$\begin{bmatrix} A \\ B_i U^{-1} \end{bmatrix} = W^{-1}W' \begin{bmatrix} A \\ B_j U^{-1} \end{bmatrix} \tag{50}$$

where $U$ is defined as in (7) and $B_i$ and $B_j$ are defined as in (13) and (14). If we post-multiply by the matrix $U$, we obtain:

$$\begin{bmatrix} A \\ B_i U^{-1} \end{bmatrix} U = W^{-1}W' \begin{bmatrix} A \\ B_j U^{-1} \end{bmatrix} U \tag{51}$$

$$\begin{bmatrix} AU \\ B_i \end{bmatrix} = W^{-1}W' \begin{bmatrix} AU \\ B_j \end{bmatrix} \tag{52}$$

Now we consider the matrix $U'$ defined as in theorem 5.2.1, and define the matrix $U'''$ as follows:

$$U'''=\begin{bmatrix} U' & O \\ O & I \end{bmatrix} \quad (53)$$

and then we pre-multiply (52) by $U'''$

$$\begin{bmatrix} U' & O \\ O & I \end{bmatrix}\begin{bmatrix} AU \\ B_i \end{bmatrix} = U'''W^{-1}W'\begin{bmatrix} AU \\ B_j \end{bmatrix} \quad (54)$$

$$\begin{bmatrix} U'AU \\ B_i \end{bmatrix} = U'''W^{-1}W'\begin{bmatrix} AU \\ B_j \end{bmatrix} \quad (55)$$

We note that the matrix $\begin{bmatrix} U'AU \\ B_i \end{bmatrix}$ is in LHNF

Now from (55), multiplying by the inverses, we obtain:

$$W'^{-1}WU'''^{-1}\begin{bmatrix} U'AU \\ B_i \end{bmatrix} = \begin{bmatrix} AU \\ B_j \end{bmatrix} \quad (56)$$

we pre-multiply again by $U'''$:

$$U'''W'^{-1}WU'''^{-1}\begin{bmatrix} U'AU \\ B_i \end{bmatrix} = \begin{bmatrix} U' & O \\ O & I \end{bmatrix}\begin{bmatrix} AU \\ B_j \end{bmatrix} \quad (57)$$

hence we obtain:

$$U'''W'^{-1}WU'''^{-1}\begin{bmatrix} U'AU \\ B_i \end{bmatrix} = \begin{bmatrix} U'AU \\ B_j \end{bmatrix} \quad (58)$$

where both the matrices $\begin{bmatrix} U'AU \\ B_j \end{bmatrix}$ and $\begin{bmatrix} U'AU \\ B_i \end{bmatrix}$ are in LHNF and the product $U'''W'^{-1}WU'''^{-1}$ is an unimodular matrix. This means that the matrix $\begin{bmatrix} U'AU \\ B_j \end{bmatrix}$ is L-A-equivalent to $\begin{bmatrix} U'AU \\ B_i \end{bmatrix}$, and therefore, by theorem 2.2.4, they must be equal, so $B_i = B_j$ and consequently $S_i = S_j$. This means that the solution $X$ belongs to one and only one of the classes $C_i$ defined in 5.2.3, therefore the set of $C_i$ is a partition of S. □

### 5.2.4 *A Formula to Express all the Solutions*

Let's suppose that $A$ is given in the equation (1), then we can express all the solutions of (1) in a way that allows us to put in evidence the free parameters and reduce the number of unknowns. The formula that we obtain reduces the initial problem to the problem of finding the values of a unimodular matrix in which the number of unknowns is lower than those of the equation (1). We will prove the following:

**Theorem 5.2.4.1** If $X_i$ and $X_j$ are two L-$A$-equivalent solutions of $\begin{vmatrix} A \\ X \end{vmatrix} = \pm d$, then the unimodular matrix $W$ such that $\begin{bmatrix} A \\ X_i \end{bmatrix} = W \begin{bmatrix} A \\ X_j \end{bmatrix}$, is in blocks lower triangular form: $\begin{bmatrix} I & O \\ \Lambda & W' \end{bmatrix}$ in which $I$ is an unit matrix, $\Lambda$ is a matrix with integer entries, $O$ is a matrix with all entries equal to 0, and $W'$ is an unimodular matrix.

*Proof.* We recall first some facts concerning the matrices in blocks lower triangular form. Let's consider the set of the real square non singular matrices in block lower triangular form, in which there are two blocks of fixed dimension, that is, of the form: $\begin{bmatrix} B_1 & O \\ R & B_2 \end{bmatrix}$ (59)

1. The determinant of a matrix in this form is given by the product of the block determinants [10, 11]

2. The product of two matrices of this form is of the same form, in fact:

$$\begin{bmatrix} B_1 & O \\ R & B_2 \end{bmatrix} \cdot \begin{bmatrix} B'_1 & O \\ R' & B'_2 \end{bmatrix} = \begin{bmatrix} B_1 B'_1 & B_1 O + O B'_2 \\ R B'_1 + B_2 R' & R O + B_2 B'_2 \end{bmatrix} = \begin{bmatrix} B_1 B'_1 & O \\ R'' & B_2 B'_2 \end{bmatrix} \quad (60)$$

3. The inverse of a matrix of this form is of the same form (we consider values in the real field), in fact given a non singular matrix, is known that the inverse is unique, and the following product shows that the inverse has the same form (Note that from 1 it follows that none of the blocks is singular):

$$\begin{bmatrix} B_1 & O \\ R & B_2 \end{bmatrix} \cdot \begin{bmatrix} B_1^{-1} & O \\ -B_2^{-1} R B_1^{-1} & B_2^{-1} \end{bmatrix} = \begin{bmatrix} I & O \\ R B_1^{-1} - B_2 B_2^{-1} R B_1^{-1} & I \end{bmatrix} = \begin{bmatrix} I & O \\ O & I \end{bmatrix} \quad (61)$$

(Since the identity matrix has the same form, it follows that the set of this kind of matrix is a group)

Be $U$ the unimodular matrix defined in (7), if we consider the two matrices:

$$M_i = \begin{bmatrix} A \\ X_i \end{bmatrix} U = \begin{bmatrix} AU \\ X_i U \end{bmatrix} = \begin{bmatrix} A' & O \\ X_{i1} & X_{i2} \end{bmatrix} \quad (62)$$

$$M_j = \begin{bmatrix} A \\ X_j \end{bmatrix} U = \begin{bmatrix} AU \\ X_j U \end{bmatrix} = \begin{bmatrix} A' & O \\ X_{j1} & X_{j2} \end{bmatrix}$$

The two matrices are clearly in blocks lower triangular form, then we can consider the equation

$$\begin{bmatrix} A \\ X_i \end{bmatrix} = W \begin{bmatrix} A \\ X_j \end{bmatrix} \quad (63)$$

and, by multiplying by $U$

$$\begin{bmatrix} A \\ X_i \end{bmatrix} U = W \begin{bmatrix} A \\ X_j \end{bmatrix} U \tag{64}$$

$$M_i = W M_j \Rightarrow M_i M_j^{-1} = W \tag{65}$$

therefore, by using the previous facts 2 and 3, the matrix $W$ is in blocks lower triangular form (but its entries are integers!). Hence we can write:

$$\begin{bmatrix} A \\ X_i \end{bmatrix} = \begin{bmatrix} W_1 & O \\ W_2 & W_3 \end{bmatrix} \begin{bmatrix} A \\ X_j \end{bmatrix} = \begin{bmatrix} W_1 A \\ W_2 A + W_3 X_j \end{bmatrix} \tag{66}$$

It means, since the matrix $A$ is full rank, that the matrix $W_1$ is the identity $I$, and the determinant of $W_3$ is 1 or -1, that is, $W_3$ is unimodular. □

If we consider that every solution $X$ is L-$A$-equivalent to one of the $S_i$, we can write

$$\begin{bmatrix} A \\ X \end{bmatrix} = \begin{bmatrix} I & O \\ \Lambda & W' \end{bmatrix} \begin{bmatrix} A \\ S_i \end{bmatrix} = \begin{bmatrix} IA + OS_i \\ \Lambda A + W' S_i \end{bmatrix} = \begin{bmatrix} A \\ \Lambda A + W' S_i \end{bmatrix} \tag{67}$$

therefore all the solutions of (1) can be expressed by the formula:

$$X = \Lambda A + W' S_i \tag{68}$$

The matrix $W$ is unique, since the matrices $X_i$ and $X_j$ are full rank. Therefore every possible matrix $\Lambda$ and every unimodular matrix $W'$ give a solution of (1), and for each solution we can find a unique $S_i$, a unique unimodular matrix $W'$ and a unique $\Lambda$ which, by formula (68), give the solution. From formula (68), and by denoting with $r$ the number of rows of the matrix $A$, and with $n$ the dimension of the matrix $\begin{bmatrix} A \\ X \end{bmatrix}$, we can see that the number of new unknowns (the entries of the unimodular matrix W') is given by $(n-r) \cdot (n-r)$, and the number of free parameters (the entries of the integer matrix $\Lambda$) is given by $r \cdot n$. The number of unknowns in (1) is given by $n(n-r)$, so the number of new unknowns in (68) is lower, since $(n-r) \cdot (n-r) < n(n-r)$.

If the matrix $A$ is not given, then the equation (1) becomes $|X| = \pm d$, the formula (68) becomes $X = W' S_i$, where the $S_i$ are the solutions already defined in theorem 5.1.1. In this case the number of unknowns is the same, since $r = 0$.

5.3 *An Example*

Let's suppose that the equation to solve is $\begin{vmatrix} 2 & 2 & -3 & 4 \\ 2 & 2 & 1 & 2 \\ x_{31} & x_{32} & x_{33} & x_{34} \\ x_{41} & x_{42} & x_{43} & x_{44} \end{vmatrix} = \pm 4$ (69)

In the example 4.2 we have seen that the greatest divisor of $A$ is 2, it divides 4, hence the equation is solvable. To find a set of solutions, we have to find the matrix $U$ defined in (7), which we have already found in example 2.4.1 :

$$U = \begin{bmatrix} 3 & 2 & -5 & -1 \\ 0 & 0 & 0 & 1 \\ -1 & 0 & 2 & 0 \\ -2 & -1 & 4 & 0 \end{bmatrix}$$

then we have to compute the inverse matrix

$$U^{-1} = \begin{bmatrix} 2 & 2 & -3 & 4 \\ 0 & 0 & 2 & -1 \\ 1 & 1 & -1 & 2 \\ 0 & 1 & 0 & 0 \end{bmatrix} \quad (70)$$

now we consider all the matrices $B_i$, defined in (13) and satisfying conditions (12) and (14). Then, considering that

$$AU = \begin{bmatrix} 2 & 2 & -3 & 4 \\ 2 & 2 & 1 & 2 \end{bmatrix} \cdot \begin{bmatrix} 3 & 2 & -5 & -1 \\ 0 & 0 & 0 & 1 \\ -1 & 0 & 2 & 0 \\ -2 & -1 & 4 & 0 \end{bmatrix} = \begin{bmatrix} 1 & 0 & 0 & 0 \\ 1 & 2 & 0 & 0 \end{bmatrix} \text{ and } k = \frac{d}{|A|} = \frac{4}{2} = 2 \quad (71)$$

all the possibile matrices $B_i$ are the following:

$$B_1 = \begin{bmatrix} 0 & 0 & 1 & 0 \\ 0 & 0 & 0 & 2 \end{bmatrix} ; B_2 = \begin{bmatrix} 0 & 0 & 2 & 0 \\ 0 & 0 & 0 & 1 \end{bmatrix} ; B_3 = \begin{bmatrix} 0 & 0 & 2 & 0 \\ 0 & 0 & 1 & 1 \end{bmatrix} ;$$

$$B_4 = \begin{bmatrix} 0 & 1 & 1 & 0 \\ 0 & 0 & 0 & 2 \end{bmatrix} ; B_5 = \begin{bmatrix} 0 & 1 & 2 & 0 \\ 0 & 0 & 0 & 1 \end{bmatrix} ; B_6 = \begin{bmatrix} 0 & 1 & 2 & 0 \\ 0 & 0 & 1 & 1 \end{bmatrix} \quad (72)$$

$$B_7 = \begin{bmatrix} 0 & 1 & 1 & 0 \\ 0 & 1 & 0 & 2 \end{bmatrix} ; B_8 = \begin{bmatrix} 0 & 1 & 2 & 0 \\ 0 & 1 & 0 & 1 \end{bmatrix} ; B_9 = \begin{bmatrix} 0 & 1 & 2 & 0 \\ 0 & 1 & 1 & 1 \end{bmatrix}$$

$$B_{10} = \begin{bmatrix} 0 & 0 & 1 & 0 \\ 0 & 1 & 0 & 2 \end{bmatrix} ; B_{11} = \begin{bmatrix} 0 & 0 & 2 & 0 \\ 0 & 1 & 0 & 1 \end{bmatrix} ; B_{12} = \begin{bmatrix} 0 & 0 & 2 & 0 \\ 0 & 1 & 1 & 1 \end{bmatrix}$$

then we can find the set of fundamental solutions by applying the formula (19)

$$S_i = B_i U^{-1}$$

$$S_1 = \begin{bmatrix} 0 & 0 & 1 & 0 \\ 0 & 0 & 0 & 2 \end{bmatrix} \begin{bmatrix} 2 & 2 & -3 & 4 \\ 0 & 0 & 2 & -1 \\ 1 & 1 & -1 & 2 \\ 0 & 1 & 0 & 0 \end{bmatrix} = \begin{bmatrix} 1 & 1 & -1 & 2 \\ 0 & 2 & 0 & 0 \end{bmatrix} \tag{73}$$

we can verify that $\begin{vmatrix} 2 & 2 & -3 & 4 \\ 2 & 2 & 1 & 2 \\ 1 & 1 & -1 & 2 \\ 0 & 2 & 0 & 0 \end{vmatrix} = 4$

By operating similarly with all the $B_i$, we obtain all the fundamental solutions:

$$S_1 = \begin{bmatrix} 1 & 1 & -1 & 2 \\ 0 & 2 & 0 & 0 \end{bmatrix} \quad S_2 = \begin{bmatrix} 2 & 2 & -2 & 4 \\ 0 & 1 & 0 & 0 \end{bmatrix} \quad S_3 = \begin{bmatrix} 2 & 2 & -2 & 4 \\ 1 & 2 & -1 & 2 \end{bmatrix}$$

$$S_4 = \begin{bmatrix} 1 & 1 & 1 & 1 \\ 0 & 2 & 0 & 0 \end{bmatrix} \quad S_5 = \begin{bmatrix} 2 & 2 & 0 & 3 \\ 0 & 1 & 0 & 0 \end{bmatrix} \quad S_6 = \begin{bmatrix} 2 & 2 & 0 & 3 \\ 1 & 2 & -1 & 2 \end{bmatrix}$$

$$S_7 = \begin{bmatrix} 1 & 1 & 1 & 1 \\ 0 & 2 & 2 & -1 \end{bmatrix} \quad S_8 = \begin{bmatrix} 2 & 2 & 0 & 3 \\ 0 & 1 & 2 & -1 \end{bmatrix} \quad S_9 = \begin{bmatrix} 2 & 2 & 0 & 3 \\ 1 & 2 & 1 & 1 \end{bmatrix} \tag{74}$$

$$S_{10} = \begin{bmatrix} 1 & 1 & -1 & 2 \\ 0 & 2 & 2 & -1 \end{bmatrix} \quad S_{11} = \begin{bmatrix} 2 & 2 & -2 & 4 \\ 0 & 1 & 2 & -1 \end{bmatrix} \quad S_{12} = \begin{bmatrix} 2 & 2 & -2 & 4 \\ 1 & 2 & 1 & 1 \end{bmatrix}$$

For each of the $S_i$ a class $C_i$ of L-$A$-equivalence is defined, and any solution $X$ can be obtained, by the formula (68):

$$\begin{bmatrix} x_{31} & x_{32} & x_{33} & x_{34} \\ x_{41} & x_{42} & x_{43} & x_{44} \end{bmatrix} = \begin{bmatrix} \lambda_{11} & \lambda_{12} \\ \lambda_{21} & \lambda_{22} \end{bmatrix} \cdot \begin{bmatrix} 2 & 2 & -3 & 4 \\ 2 & 2 & 1 & 2 \end{bmatrix} + \begin{bmatrix} w'_{11} & w'_{12} \\ w'_{21} & w'_{22} \end{bmatrix} \cdot [S_i] \tag{75}$$

where $\begin{vmatrix} w'_{11} & w'_{12} \\ w'_{21} & w'_{22} \end{vmatrix} = \pm 1$, $S_i$ takes value in the set (74) and the $\lambda_{ij}$ can take any integer value.

5.4 *The Number of L-A-Equivalence Classes*

We take into consideration the total amount of the different $S_i$, the fundamental solutions given $A$, in §5.4. Let's suppose that the equation $\begin{vmatrix} A \\ X \end{vmatrix} = \pm d$ is solvable, i.e $|A|$ divides $d$, and set $d = p_1^{n_1} p_2^{n_2} .. p_l^{n_l} |A|$, in which the $p_i$ are different prime numbers, and $k$ the row number of the matrix $X$. We will prove the following:

**Theorem 5.4** The number N of L-A-equivalence classes (or equivalently the cardinality of the fundamental solutions set), is given by:

$$N = |A|^k (\sum_{j_0=1}^{k}\sum_{j_1=j_0}^{k}\sum_{j_2=j_1}^{k} \cdots \sum_{j_{n_1}=j_{n_1-1}}^{k} p_1^{n_1 k - j_0 - j_1 \cdots - j_{n_1}})(\sum_{j_0=1}^{k}\sum_{j_1=j_0}^{k}\sum_{j_2=j_1}^{k} \cdots \sum_{j_{n_2}=j_{n_2-1}}^{k} p_2^{n_2 k - j_0 - j_1 \cdots - j_{n_2}}) \cdots (\sum_{j_0=1}^{k}\sum_{j_1=j_0}^{k}\sum_{j_2=j_1}^{k} \cdots \sum_{j_{n_l}=j_{n_l-1}}^{k} p_l^{n_l k - j_0 - j_1 \cdots - j_{n_l}})$$

(76)

in which $k$ is the number of rows of $X$

*Proof.*

Let's start proving the factor $|A|^k$:

if we take into consideration the hermitian form of $\begin{bmatrix} A \\ X \end{bmatrix}$ :
$\begin{bmatrix} a'_{11} & 0 & .. & .. & .. & 0 \\ .. & .. & 0 & .. & .. & 0 \\ a'_{r1} & & a'_{rr} & 0 & .. & 0 \\ b_{l1} & b_{l2} & .. & b_{ll} & .. & 0 \\ .. & .. & & & .. & 0 \\ b_{n1} & & & & & b_{nn} \end{bmatrix}$

every $b_{ij}$, with $j \leq r$, can take $a_{ii}$ possible integer values from 0 to $a_{ii} - 1$. Since in every column there are $k$ possible entries, all the possibilities are given by $a_{11}^k a_{22}^k \ldots a_{rr}^k = (a_{11} a_{22} \ldots a_{rr})^k = |A|^k$.

Now if we consider the entries $b_{ii}$, with $r < i \leq n$, we have: $\prod_{i=r}^{n} b_{ii} = \frac{d}{|A|} = p_1^{n_1} p_2^{n_2} \ldots p_l^{n_l}$

If we consider the square matrix B (lower triangular) with the diagonal entries $b_{ii}$, all the possible integer values of the column i are given by $(b_{ii})^{k-i}$ and considering all the diagonal entries $b_{ii}$ the number N of possible matrices will be the product $N = (b_{11})^{k-1} \cdot (b_{22})^{k-2} \cdots (b_{kk})^0$. If we take into consideration the prime factorization can write:

$$N = (p_1^{e_{11}})^{k-1} \cdot (p_2^{e_{12}})^{k-1} \cdots (p_1^{e_{21}})^{k-2} \cdot (p_2^{e_{22}})^{k-2} \cdots (p_1^{e_{31}})^{k-3} \cdot (p_2^{e_{32}})^{k-3} \cdots (p_1^{e_{k1}})^0 \cdot (p_2^{e_{k2}})^0$$

and by reordering the terms:

$$N = (p_1^{e_{11}})^{k-1} \cdot (p_1^{e_{21}})^{k-2} \cdot (p_1^{e_{31}})^{k-3} \cdots (p_2^{e_{12}})^{k-1} \cdot (p_2^{e_{22}})^{k-2} \cdot (p_2^{e_{32}})^{k-3} \cdots (p_l^{e_{l1}})^{k-1} \cdot (p_l^{e_{l2}})^{k-2} \cdot (p_l^{e_{l3}})^{k-3} \cdots$$

It means that in order to find the total number of the possible matrices B we can find for each factor $p_i^{n_i}$ the number of possible square lower triangular matrices of order k, with determinant equal to $p_i^{n_i}$ and then we multiply all these numbers.

Now, if we consider the factor $p_i^{n_i}$, and consider a lower triangular matrix of order k and determinant $p_i^{n_i}$, for each $p_i$ set in the diagonal on the column (or row) number j, we have a contribution (by multiplication) of $p_i^{k-j}$, so, by considering all possible disposition on the diagonal, we obtain that the total number of such matrices is given by $\sum_{j_0=1}^{k}\sum_{j_1=j_0}^{k}\sum_{j_2=j_1}^{k}..\sum_{j_{n_i}=j_{n_i-1}}^{k} p_i^{k-j_0} p_i^{k-j_1}..p_i^{k-j_{n_i}}$, which can be written as: $\sum_{j_0=1}^{k}\sum_{j_1=j_0}^{k}\sum_{j_2=j_1}^{k}..\sum_{j_{n_i}=j_{n_i-1}}^{k} p_i^{n_i k - j_0 - j_1 - j_2 ... - j_{n_i}}$

Finally, by considering all the factors $p_i^{n_i}$, we can write the formula which gives the number of all the possible L-A-equivalence classes:

$$N = |A|^k (\sum_{j_0=1}^{k}\sum_{j_1=j_0}^{k}\sum_{j_2=j_1}^{k}..\sum_{j_{n_1}=j_{n_1-1}}^{k} p_1^{n_1 k - j_0 - j_1 .. - j_{n_1}})(\sum_{j_0=1}^{k}\sum_{j_1=j_0}^{k}\sum_{j_2=j_1}^{k}..\sum_{j_{n_2}=j_{n_2-1}}^{k} p_2^{n_2 k - j_0 - j_1 .. - j_{n_2}})..(\sum_{j_0=1}^{k}\sum_{j_1=j_0}^{k}\sum_{j_2=j_1}^{k}..\sum_{j_{nl}=j_{nl-1}}^{k} p_l^{n_l k - j_0 - j_1 .. - j_{nl}})$$

□

5.4.1 *Example*

If we consider the equation (69), we have that, $|A|=2$, $k=2$

clearly we have $p_1 = 2$ and $n_1 = 1$. By applying the formula (76) we have:

$$N = |A|^k (\sum_{j=1}^{2} 2^{2-j}) = 2^2 (2^1 + 2^0) = 4 \cdot 3 = 12$$

## 6. Generation of Unimodular Matrices

If we consider the equation (68), it's clear that in order to find explicitly all the solutions we need a way to build all the possible unimodular matrices. Some known results on the $GL(n,\mathbb{Z})$ group of $n \times n$ integral matrices can be useful. It's known that 4 matrices suffice to generate this group [14].

The number of generators can be reduced to the following 2, as Trott has proven [15]:

$$U_1 = \begin{bmatrix} 1 & 0 & 0 & .. & .. & 0 \\ 1 & 1 & 0 & .. & .. & 0 \\ 0 & 0 & 1 & 0 & .. & 0 \\ .. & .. & 0 & .. & .. & 0 \\ .. & .. & .. & .. & 1 & 0 \\ 0 & 0 & 0 & .. & 0 & 1 \end{bmatrix} \quad U = \begin{bmatrix} 0 & 1 & 0 & .. & .. & 0 \\ 0 & 0 & 1 & .. & .. & 0 \\ 0 & 0 & 0 & 1 & .. & 0 \\ .. & & .. & 0 & .. & .. & 0 \\ .. & & .. & .. & .. & 0 & 1 \\ (-1)^n & 0 & 0 & .. & 0 & 0 \end{bmatrix}$$

Hence, going back to the formula (68), that gives the general solutions of equation $\begin{vmatrix} A \\ X \end{vmatrix} = \pm d$,

$$X = \Lambda A + W'S_i$$

we can write any unimodular matrix W' as product of U and $U_1$, $W' = \prod_{i=1}^{k}(U^{a_i} U_1^{b_i})$

where $k$ is natural and $a_i$, $b_i$ are integers. Therefore the general solution can be written as:

$$X = \Lambda A + \left( \prod_{j=1}^{k} U^{a_j} U_1^{b_j} \right) S_i \qquad (77)$$

6.1 *Example*

We can consider the expression (75):

$$\begin{bmatrix} x_{31} & x_{32} & x_{33} & x_{34} \\ x_{41} & x_{42} & x_{43} & x_{44} \end{bmatrix} = \begin{bmatrix} \lambda_{11} & \lambda_{12} \\ \lambda_{21} & \lambda_{22} \end{bmatrix} \cdot \begin{bmatrix} 2 & 2 & -3 & 4 \\ 2 & 2 & 1 & 2 \end{bmatrix} + \begin{bmatrix} w'_{11} & w'_{12} \\ w'_{21} & w'_{22} \end{bmatrix} \cdot [S_i]$$

where $\begin{vmatrix} w'_{11} & w'_{12} \\ w'_{21} & w'_{22} \end{vmatrix} = \pm 1$ and write the unimodular matrix as product of the generators. In this case the generators are: $U = \begin{bmatrix} 0 & 1 \\ 1 & 0 \end{bmatrix}$ and $U_1 = \begin{bmatrix} 1 & 0 \\ 1 & 1 \end{bmatrix}$

therefore we can obtain all solutions of equation (69) by the formula (77):

$$\begin{bmatrix} x_{31} & x_{32} & x_{33} & x_{34} \\ x_{41} & x_{42} & x_{43} & x_{44} \end{bmatrix} = \begin{bmatrix} \lambda_{11} & \lambda_{12} \\ \lambda_{21} & \lambda_{22} \end{bmatrix} \cdot \begin{bmatrix} 2 & 2 & -3 & 4 \\ 2 & 2 & 1 & 2 \end{bmatrix} + \left( \prod_{j=1}^{k} \begin{bmatrix} 0 & 1 \\ 1 & 0 \end{bmatrix}^{a_j} \cdot \begin{bmatrix} 1 & 0 \\ 1 & 1 \end{bmatrix}^{b_j} \right) \cdot [S_i]$$

where the fundamental solutions $S_i$ takes values in the set (74), k is a natural number $\lambda_{ij}$, $a_j$, $b_j$ are integers.

6.2 *Another Example*

Let's solve the equation : $\begin{vmatrix} 1 & 2 & 2 & 0 & 0 \\ -1 & 1 & 3 & 0 & 1 \\ x_{31} & x_{32} & x_{33} & x_{34} & x_{35} \\ x_{41} & x_{42} & x_{43} & x_{44} & x_{45} \\ x_{51} & x_{52} & x_{53} & x_{54} & x_{55} \end{vmatrix} = \pm 4$

1- First we compute the gcd to establish if it's solvable: gcd (3, 5, 0, 1, 4, 0, 2, 0, 2, 0)=1, 1 divides 4, hence the equation is solvable.

2- Then we find the RHNF of $A$, the matrix $U$ and $U^{-1}$:

$$\begin{bmatrix} 1 & 2 & 2 & 0 & 0 \\ -1 & 1 & 3 & 0 & 1 \end{bmatrix} \cdot \begin{bmatrix} -1 & -2 & 8 & 0 & 2 \\ 1 & 1 & -5 & 0 & -1 \\ 0 & 0 & 1 & 0 & 0 \\ 0 & 0 & 0 & 1 & 0 \\ -2 & -2 & 10 & 0 & 3 \end{bmatrix} = \begin{bmatrix} 1 & 0 & 0 & 0 & 0 \\ 0 & 1 & 0 & 0 & 0 \end{bmatrix}$$

$$U^{-1} = \begin{bmatrix} -1 & -2 & 8 & 0 & 2 \\ 1 & 1 & -5 & 0 & -1 \\ 0 & 0 & 1 & 0 & 0 \\ 0 & 0 & 0 & 1 & 0 \\ -2 & -2 & 10 & 0 & 3 \end{bmatrix}^{-1} = \begin{bmatrix} 1 & 2 & 2 & 0 & 0 \\ -1 & 1 & 3 & 0 & 1 \\ 0 & 0 & 1 & 0 & 0 \\ 0 & 0 & 0 & 1 & 0 \\ 0 & 2 & 0 & 0 & 1 \end{bmatrix}$$

3- Find the fundamental solutions

The cardinality of the fundamental solutions set is given by
$$N = 1^3 (\sum_{j_0=1}^{3} \sum_{j_1=j_0}^{3} 2^{6-j_0-j_1}) = (2^{6-1-1} + 2^{6-1-2} + 2^{6-1-3}) + (2^{6-2-2} + 2^{6-2-3}) + (2^{6-3-3}) = (16+8+4)+(4+2)+1 = 35$$

The 35 $B_i$ matrices are:

$$\begin{bmatrix} 0 & 0 & 1 & 0 & 0 \\ 0 & 0 & 0 & 1 & 0 \\ 0 & 0 & 0 & 0 & 4 \end{bmatrix} \begin{bmatrix} 0 & 0 & 1 & 0 & 0 \\ 0 & 0 & 0 & 4 & 0 \\ 0 & 0 & 0 & 0 & 1 \end{bmatrix} \begin{bmatrix} 0 & 0 & 1 & 0 & 0 \\ 0 & 0 & 0 & 4 & 0 \\ 0 & 0 & 0 & 1 & 1 \end{bmatrix} \begin{bmatrix} 0 & 0 & 1 & 0 & 0 \\ 0 & 0 & 0 & 4 & 0 \\ 0 & 0 & 0 & 2 & 1 \end{bmatrix} \begin{bmatrix} 0 & 0 & 1 & 0 & 0 \\ 0 & 0 & 0 & 4 & 0 \\ 0 & 0 & 0 & 3 & 1 \end{bmatrix}$$

$$\begin{bmatrix} 0 & 0 & 4 & 0 & 0 \\ 0 & 0 & 0 & 1 & 0 \\ 0 & 0 & 0 & 0 & 1 \end{bmatrix} \begin{bmatrix} 0 & 0 & 4 & 0 & 0 \\ 0 & 0 & 1 & 1 & 0 \\ 0 & 0 & 0 & 0 & 1 \end{bmatrix} \begin{bmatrix} 0 & 0 & 4 & 0 & 0 \\ 0 & 0 & 2 & 1 & 0 \\ 0 & 0 & 0 & 0 & 1 \end{bmatrix} \begin{bmatrix} 0 & 0 & 4 & 0 & 0 \\ 0 & 0 & 3 & 1 & 0 \\ 0 & 0 & 0 & 0 & 1 \end{bmatrix} \begin{bmatrix} 0 & 0 & 4 & 0 & 0 \\ 0 & 0 & 0 & 1 & 0 \\ 0 & 0 & 1 & 0 & 1 \end{bmatrix}$$

$$\begin{bmatrix} 0 & 0 & 4 & 0 & 0 \\ 0 & 0 & 1 & 1 & 0 \\ 0 & 0 & 1 & 0 & 1 \end{bmatrix} \begin{bmatrix} 0 & 0 & 4 & 0 & 0 \\ 0 & 0 & 2 & 1 & 0 \\ 0 & 0 & 1 & 0 & 1 \end{bmatrix} \begin{bmatrix} 0 & 0 & 4 & 0 & 0 \\ 0 & 0 & 3 & 1 & 0 \\ 0 & 0 & 1 & 0 & 1 \end{bmatrix} \begin{bmatrix} 0 & 0 & 4 & 0 & 0 \\ 0 & 0 & 0 & 1 & 0 \\ 0 & 0 & 2 & 0 & 1 \end{bmatrix} \begin{bmatrix} 0 & 0 & 4 & 0 & 0 \\ 0 & 0 & 1 & 1 & 0 \\ 0 & 0 & 2 & 0 & 1 \end{bmatrix}$$

$$\begin{bmatrix} 0 & 0 & 4 & 0 & 0 \\ 0 & 0 & 2 & 1 & 0 \\ 0 & 0 & 2 & 0 & 1 \end{bmatrix} \begin{bmatrix} 0 & 0 & 4 & 0 & 0 \\ 0 & 0 & 3 & 1 & 0 \\ 0 & 0 & 2 & 0 & 1 \end{bmatrix} \begin{bmatrix} 0 & 0 & 4 & 0 & 0 \\ 0 & 0 & 0 & 1 & 0 \\ 0 & 0 & 3 & 0 & 1 \end{bmatrix} \begin{bmatrix} 0 & 0 & 4 & 0 & 0 \\ 0 & 0 & 1 & 1 & 0 \\ 0 & 0 & 3 & 0 & 1 \end{bmatrix} \begin{bmatrix} 0 & 0 & 4 & 0 & 0 \\ 0 & 0 & 2 & 1 & 0 \\ 0 & 0 & 3 & 0 & 1 \end{bmatrix}$$

$$\begin{bmatrix} 0 & 0 & 4 & 0 & 0 \\ 0 & 0 & 3 & 1 & 0 \\ 0 & 0 & 3 & 0 & 1 \end{bmatrix} \begin{bmatrix} 0 & 0 & 1 & 0 & 0 \\ 0 & 0 & 0 & 2 & 0 \\ 0 & 0 & 0 & 0 & 2 \end{bmatrix} \begin{bmatrix} 0 & 0 & 1 & 0 & 0 \\ 0 & 0 & 0 & 2 & 0 \\ 0 & 0 & 0 & 1 & 2 \end{bmatrix} \begin{bmatrix} 0 & 0 & 2 & 0 & 0 \\ 0 & 0 & 0 & 1 & 0 \\ 0 & 0 & 0 & 0 & 2 \end{bmatrix} \begin{bmatrix} 0 & 0 & 2 & 0 & 0 \\ 0 & 0 & 1 & 1 & 0 \\ 0 & 0 & 0 & 0 & 2 \end{bmatrix}$$

$$\begin{bmatrix} 0 & 0 & 2 & 0 & 0 \\ 0 & 0 & 1 & 1 & 0 \\ 0 & 0 & 1 & 0 & 2 \end{bmatrix} \begin{bmatrix} 0 & 0 & 2 & 0 & 0 \\ 0 & 0 & 0 & 1 & 0 \\ 0 & 0 & 1 & 0 & 2 \end{bmatrix} \begin{bmatrix} 0 & 0 & 2 & 0 & 0 \\ 0 & 0 & 0 & 2 & 0 \\ 0 & 0 & 0 & 0 & 1 \end{bmatrix} \begin{bmatrix} 0 & 0 & 2 & 0 & 0 \\ 0 & 0 & 0 & 2 & 0 \\ 0 & 0 & 0 & 1 & 1 \end{bmatrix} \begin{bmatrix} 0 & 0 & 2 & 0 & 0 \\ 0 & 0 & 1 & 2 & 0 \\ 0 & 0 & 0 & 0 & 1 \end{bmatrix}$$

$$\begin{bmatrix} 0 & 0 & 2 & 0 & 0 \\ 0 & 0 & 1 & 2 & 0 \\ 0 & 0 & 0 & 1 & 1 \end{bmatrix} \begin{bmatrix} 0 & 0 & 2 & 0 & 0 \\ 0 & 0 & 0 & 2 & 0 \\ 0 & 0 & 1 & 0 & 1 \end{bmatrix} \begin{bmatrix} 0 & 0 & 2 & 0 & 0 \\ 0 & 0 & 0 & 2 & 0 \\ 0 & 0 & 1 & 1 & 1 \end{bmatrix} \begin{bmatrix} 0 & 0 & 2 & 0 & 0 \\ 0 & 0 & 1 & 2 & 0 \\ 0 & 0 & 1 & 0 & 1 \end{bmatrix} \begin{bmatrix} 0 & 0 & 2 & 0 & 0 \\ 0 & 0 & 1 & 2 & 0 \\ 0 & 0 & 1 & 1 & 1 \end{bmatrix}$$

therefore the set of fundamental solutions is given by $S_i = B_i \, U^{-1}$

4- Write all the solutions :

$$\begin{bmatrix} x_{31} & x_{32} & x_{33} & x_{34} & x_{35} \\ x_{41} & x_{42} & x_{43} & x_{44} & x_{45} \\ x_{51} & x_{52} & x_{53} & x_{54} & x_{55} \end{bmatrix} =$$

$$= \begin{bmatrix} \lambda_{11} & \lambda_{12} \\ \lambda_{21} & \lambda_{22} \\ \lambda_{31} & \lambda_{32} \end{bmatrix} \cdot \begin{bmatrix} 1 & 2 & 2 & 0 & 0 \\ -1 & 1 & 3 & 0 & 1 \end{bmatrix} + \left( \prod_{j=1}^{k} \begin{bmatrix} 1 & 0 & 0 \\ 1 & 1 & 0 \\ 0 & 0 & 1 \end{bmatrix}^{a_j} \cdot \begin{bmatrix} 0 & 1 & 0 \\ 0 & 0 & 1 \\ -1 & 0 & 0 \end{bmatrix}^{b_j} \right) \cdot B_i \cdot \begin{bmatrix} 1 & 2 & 2 & 0 & 0 \\ -1 & 1 & 3 & 0 & 1 \\ 0 & 0 & 1 & 0 & 0 \\ 0 & 0 & 0 & 1 & 1 \\ 0 & 2 & 0 & 0 & 1 \end{bmatrix}$$

where $k$ is a natural number $\lambda_{ij}$ $a_j$, $b_j$ are integer.

## 7. Concluding Remarks

In this paper we have considered the diophantine equation $\left| \begin{matrix} A \\ X \end{matrix} \right| = \pm d$ where $A$ is full rank and $d$ is not zero. The equation is solvable if and only if the gcd of the matrix $A$ divides $d$. If the equation is solvable we have shown how to build a set of solutions $S_i$, which we have defined fundamental solutions. Then we have defined an equivalence relation which gives rise to a partition of the set of all solutions in which each $S_i$ belongs to a different class. We have given a formula to express all the solutions in which we two generators $U$ and $U_1$ of GL(n, Z) are used:

$$X = \Lambda A + \left(\prod_{j=1}^{k} U^{a_j} U_1^{b_j}\right) S_i$$

Finally we have given a formula to compute the number of the equivalence classes. In this paper we considered $d \neq 0$ and much of what we proved is valid only if the matrix $\begin{bmatrix} A \\ X \end{bmatrix}$ is not singular.